\documentclass[reqno,11pt]{amsart}
\usepackage{amsmath}
\usepackage{graphicx, color, enumerate}
\newcommand*{\rom}[1]{\expandafter\@slowromancap\romannumeral #1@}
\usepackage{amscd}
\usepackage{amsmath}
\usepackage{amsfonts}
\usepackage{amssymb}
\usepackage{mathrsfs}
\usepackage{latexsym}

\newcommand{\be}{\begin{equation}}
\newcommand{\ee}{\end{equation}}

\newenvironment{pfn}{\noindent{\it Proof. }}{\rule{2mm}{2mm}\medskip}

\newcommand{\R}{\mathbb{R}}

\newcommand{\E}{\mathbb{E}}

\newcommand{\h}{\mathbb{H}}
\newcommand{\dyle}{\displaystyle}
\renewcommand{\b }{\beta }

\newcommand{\D }{\Delta }

\renewcommand{\l }{\lambda }
\renewcommand{\L }{\Lambda }

\newcommand{\bra}{\langle} %% apre prod. scal. in Rn
\newcommand{\ket}{\rangle} %%chiude prod. scal.

\newcommand{\cM}{{\mathcal{M}}}

\newcommand{\cN}{{\mathcal{N}}}

\newcommand{\intN}{\int_{\R^N}}
\newcommand{\intR}{\int_\R}

\newcommand{\wt}{\widetilde}

\newcommand{\bu}{{\bf u}}
\newcommand{\bv}{{\bf v}}

\newcommand{\bw}{{\bf w}}
\newcommand{\bo}{{\bf 0}}

\newcommand{\bh}{{\bf h}}

\newtheorem{Theorem}{Theorem}
\newtheorem{Corollary}[Theorem]{Corollary}
\newtheorem{Lemma}[Theorem]{Lemma}
\newtheorem{Proposition}[Theorem]{Proposition}

\newtheorem{remark}[Theorem]{Remark}
\newtheorem{remarks}[Theorem]{Remarks}
\newtheorem{example}[Theorem]{Example}
\newtheorem{examples}[Theorem]{Examples}
\newenvironment{Remark}{\begin{remark}}{\end{remark}}
\newenvironment{Remarks}{\begin{remarks}
}{\end{remarks}}

        \usepackage{latexsym}
        \usepackage{amsmath}
        \usepackage{amsfonts}
        \usepackage{amssymb}
\author{P.~\'Alvarez-Caudevilla, Eduardo Colorado and Rasiel Fabelo}

\address{Universidad Carlos III de Madrid,
Av. Universidad 30, 28911-Legan\'es, Spain \& Instituto de Ciencias Matem\'aticas,
ICMAT (CSIC-UAM- UC3M-UCM), C/Nicol\'as Cabrera 15, 28049 Madrid, Spain} \email{pacaudev@math.uc3m.es, pablo.caudevilla@icmat.es}

\address{Universidad Carlos III de Madrid,
Av. Universidad 30, 28911-Legan\'es, Spain \& Instituto de Ciencias Matem\'aticas,
ICMAT (CSIC-UAM- UC3M-UCM), C/Nicol\'as Cabrera 15, 28049 Madrid, Spain} \email{ecolorad@math.uc3m.es, eduardo.colorado@icmat.es}

\address{Universidad Carlos III de Madrid,
Av. Universidad 30, 28911-Legan\'es, Spain} \email{rfabelo@math.uc3m.es}


\date{\today}

\title[Coupled higher-order nonlinear Schr\"odinger-KdV equations]{A higher order system of some coupled nonlinear Schr\"odinger and
Korteweg-de Vries equations}

\begin{document}
{\sf \maketitle

\noindent {\small  {\bf Abstract.} We prove existence and
multiplicity of bound and ground state solutions, under appropriate
conditions on the parameters, for  a bi-harmonic stationary system
of coupled nonlinear Schr\"odinger--Korteweg-de Vries equations.

\section{Introduction}

Recently in \cite{c2,c3} has been analyzed a system of coupled
nonlinear Schr\"odinger--Korteweg-de Vries equations
\begin{equation}\label{NLS-KdV edu1}
\left\{\begin{array}{rcl}
if_t +  f_{xx} + |f|^2f+ \b fg & = & 0\\
g_t+g_{xxx}+gg_x+\tfrac12\beta(|f|^2)_x& = & 0,
\end{array}\right.
\end{equation}
with $f=f(x,t)\in \mathbb{C}$, $g=g(x,t)\in \mathbb{R}$, and
$\b\in\R$ a coupling parameter. This system appears in phenomena of
interactions between short and long dispersive waves, arising in
fluid mechanics, such as  the interactions of capillary - gravity
water waves \cite{KaSK}. Indeed, $f$ represents the short-wave,
while $g$ stands for the long-wave; see references
\cite{aa,c2,c3,cl,fo} for further details on similar system.
Moreover, the interaction between long and short waves appears in
magnetised plasma \cite{Kar75}, \cite{YaOi} and in many physical
phenomena as well, such that Bose-Einstein condensates \cite{Ch}.

The solutions studied  in papers \cite{c2,c3} (see also
\cite{dfo1,dfo2}) are taken as solitary traveling waves, i.e.,
\begin{equation}
(f(x,t),g(x,t))=(e^{i\omega
t}e^{i\frac{c}{2}x}u(x-ct),v(x-ct)),\qquad \text{ where }u,v\text{
are real functions.}
\end{equation}
Choosing  $\lambda_1=\omega+\frac{c^2}{4}$ and $\l_2=c$, then $u,\,
v$ are solutions of the following stationary system
 \be \left\lbrace
\begin{array}{ccc}\label{eq:NLS-KdV2 segundo}
-u''+\l_1u & =& u^3+\beta uv\\
-v''+\l_2v & =& \frac{1}{2}v^2+\frac{1}{2}\beta u^2.
\end{array}
 \right.
\ee

In the present  work we  analyze the existence of solutions of a
higher order system coming from \eqref{NLS-KdV edu1}. More
precisely, we consider the following system
\begin{equation}\label{NLS-KdV edu}
\left\{\begin{array}{rcl}
if_t -  f_{xxxx} + |f|^2f+ \b fg & = & 0\\
g_t-g_{xxxxx}+\tfrac{1}{2}(|g|g)_x+\tfrac12\beta(|f|^2)_x& = & 0.
\end{array}\right.
\end{equation}
Looking for ``standing-traveling"\footnote{This is the first time,
up to our knowledge, that the interaction of standing waves and
traveling waves is analyzed in the mathematical literature.} wave
solutions of the form
$$
(f(x,t),g(x,t))=\left(e^{i\l_1 t} u(x),v(x-\l_2t)\right),\qquad
\text{ where }u,v\text{ are real functions,}
$$
then we arrive at the fourth-order stationary system
 \be \left\lbrace
\begin{array}{ccc}\label{eq:NLS-KdV222}
u^{(iv} +\l_1u & =& u^3+\beta uv\\
v^{(iv}+\l_2v & =& \frac{1}{2}|v|v+\frac{1}{2}\beta u^2,
\end{array}
 \right.
\ee where $w^{(iv}$ denotes the fourth derivative of $w$. Although
system \eqref{NLS-KdV edu} has sense only in dimension $1$, passing
to the stationary system \eqref{eq:NLS-KdV222}, it makes sense to
consider it in higher dimensional cases, as the following,
\be
\left\lbrace
\begin{array}{ccc}\label{eq:NLS-KdV2}
\D^2u +\l_1u & =& u^3+\beta uv\\
\D^2v+\l_2v & =& \frac{1}{2}|v|v+\frac{1}{2}\beta u^2,
\end{array}
 \right.
\ee
 where $u,v\in W^{2,2}(\R^N)$, $1\le N\le 7$, $\l_j>0$ with
$j=1,2$ and $\beta>0$ is the coupling parameter.

Recently, other similar fourth-order systems studying the interaction of
coupled nonlinear Schr\"odinger equations have
appeared; see \cite{pev}, where the
coupling terms have the same homogeneity as the nonlinear terms. Note that, as far as we know there is not
any previous mathematical work analyzing a higher order system with the nonlinear and coupling terms considered here in the system \eqref{eq:NLS-KdV2}.

Here we first analyze the dimensional case
$2\le N\le 7$ in the radial framework (see subsection \ref{sec:4.1}) by using the compactness described in Remark
\ref{re:functionl constrain y PS condition}-$(ii)$. The one dimensional case is studied in subsection \ref{sec:4.2} where we use a
 measure Lemma due to P. L. Lions \cite{lions2} to circumvent the lack of compactness.

System \eqref{eq:NLS-KdV2} has a non-negative semi-trivial solution,
$\bv_2=(0,V_2)$ defined in Remark \ref{semi-trivial}. Then in order
to find non-negative bound or ground state solutions we need to
check that they are different from $\bv_2$. To be more precise, we
prove that there exists a positive critical value of the coupling
parameter $\beta$, denoted by $\L$ defined by
 \eqref{Lambda},
such that the associated functional constrained to the corresponding Nehari manifold possesses
 a positive global minimum, which is a critical point with energy below the energy of the semi-trivial solution under the following
 hypotheses:
 either $\beta>\Lambda$ or  $\beta>0$ and $\lambda_2\gg 1$.
Furthermore, we find a mountain pass critical point if $\beta<\L$ and $\lambda_2\gg 1$.

\

The paper is organized as follows. In Section \ref{sec:2} we
introduce the notation, establish the functional framework,
 define the Nehari manifold and study its properties.
Section \ref{sec:4} is devoted to prove the main results of the
paper. It is divided into two subsections, in the first one
(Subsection \ref{sec:4.1}) we study the high-dimensional case ($2\le
N\le 7$), while the second one  (Subsection \ref{sec:4.2}) deals
with the one-dimensional case.

\section{Functional setting, notation and Nehari manifold}\label{sec:2}
Let $E$ be the Sobolev space $W^{2,2}(\R^N)$ then, we define the following equivalent norms and scalar products:

\[\bra u,v\ket_j:=\int_{\R^N}\Delta u\cdot\Delta v\, dx+\l_j\int_{\R^N}uv\,dx,\qquad\|u\|_j^2:=\bra u,u\ket_j\qquad j=1,2.\]
Let us define the product Sovolev space $\E:=E\times E$ and denote its elements by $\bu=(u,v)$ with $\bo=(0,0)$.
We will take the inner product in $\E$ as follow,
\begin{equation}\label{inner product}
\bra\bu_1,\bu_2\ket:=\langle u_1,u_2\rangle_1+\bra v_1,v_2 \ket_2,
\end{equation}
which induces the following norm
$$
\|\bu\|:=\sqrt{\|u\|_1^2+\|v\|_2^2}.
$$
Moreover, for $\bu=(u,v)\in \E$, the notation $\bu\geq \bo$, resp. $\bu>\bo$, means that $u,v\geq 0$, resp. $u,v>0$.
We denote by $H$ the space of radially symmetric functions in $E$, and $\h:=H\times H$.
In addition, we define energy functionals associated to system \eqref{eq:NLS-KdV2} by
\begin{equation}\label{Phi22}
\Phi (\bu)= I_1(u)+I_2(v)- \tfrac 12\b \int_{\R^N} u^2v\,dx,\qquad \bu\in \E,
\end{equation}
where
$$
I_1(u)=\tfrac 12 \|u\|_1^2 -\tfrac 14\, \int_{\R^N} u^4\,dx,\qquad I_2(v)=\tfrac 12 \|v\|_2^2 -\tfrac 16\, \int_{\R^N} |v|^3\,dx,\qquad u,\, v\in E,
$$
are the energy functional associated to the  uncoupled  equations in \eqref{eq:NLS-KdV2}.

\begin{Remark}
We can easily see that the functional $\Phi$ is not bounded bellow
on $\mathbb{E}$. Thus, we are going to work on the so called Nehari
manifold which is a natural constraint for the functional $\Phi$,
and even more the functional constrained to the Nehari manifold is
bounded below.
\end{Remark}
We define
\be\label{eq:Psi}
\Psi(\bu)=\Phi'(\bu)[\bu]=\|\bu\|^2-\int_{\R^N} u^4 \,dx-\tfrac{1}{2} \int_{\R^N} |v|^3 \,dx -\tfrac{3}{2}\beta\int_{\R^N} u^2v \,dx.
\ee
Using the previous definition, the Nehari manifold is given by
\be \label{gemanif}
\mathcal{M} =\{ \bu\in \E\setminus\{\bo\}: \Psi(\bu)=0\}.
\ee

 This manifold will be used in order to deal with the one dimensional case in subsection \ref{sec:4.2}, in which there is no compactness, see Remark \ref{re:functionl constrain y PS condition}-$(ii)$.

In the dimensional case $2\le N\le 7$, we restrict the Nehari
Manifold to the radial setting, denoting it as
 \be\label{ramainf} \cN =\{
\bu\in \h\setminus\{\bo\}: \Psi(\bu)=0\}. \ee
Furthermore, differentiating expression \eqref{eq:Psi} yields
\begin{equation}\label{eq:phi'}
\Psi'(\bu)[\bu]=2\|\bu\|^2-4\int_{\R^N} u^4 \,dx-\tfrac{3}{2} \int_{\R^N} |v|^3 \,dx  -\tfrac{9}{2}\beta\int_{\R^N} u^2v \,dx .
\end{equation}
\begin{Remark}
All the properties we are going to prove in this section are satisfied for both $\cM$ and $\cN$, but the Palais-Smale  condition, in Lemma \ref{Lemma PS},
 is only satisfied for $\Phi$ on $\cN$, because of working on the radial setting, see again Remark \ref{re:functionl constrain y PS condition}-$(ii)$.
 To be short, we are going to demonstrate the following properties for the Nehari manifold $\cN$.
\end{Remark}
Using the fact that $\Psi(\bu)=0$ for any $\bu\in\cN$, we have
\begin{equation}\label{eq:delta psi restringida a N}
\Psi'(\bu)[\bu]=\Psi'(\bu)[\bu]-3\Psi(\bu)=- \|\bu \|^2-\int_{\R^N} u^4\,dx\,<0,\quad\forall\, \bu\in \cN.
\end{equation}
Then, $\cN$ is a locally smooth manifold near any point $\bu\neq  0$ with $\Psi(\bu)=0$.
Taking  the  derivative of the functional $\Phi$, we find
 \[
\Phi'(\bu)[\bh]=I_1'(u)[h_1]+I_2'(v)[h_2]-\beta\int_{\R^N} uvh_1\,dx-\tfrac{1}{2}\beta\int_{\R^N} u^2h_2 \,dx,
\]
The second derivative of $\Phi $ is given by
\[
\Phi''(\bu)[\bh]^2=\|\bh\|^2- 3\int_{\R^N} u^2h_1^2\,dx- \int_{\R^N}
|v|h_2^2\,dx-\beta\int_{\R^N} vh_1^2\,dx-2\beta\int_{\R^N}
uh_1h_2\,dx.
\]
It satisfies
$$\Phi''(\bo)[\bh]^2=\|\bh\|^2,$$
which is positive definite, so that $\bo$ is a strict minimum
critical point for $\Phi$. As a consequence, we have that $\cN$ is a
smooth complete manifold, and there exists a constant $\rho>0$ such
that \be\label{eq:away from zero}
\|\bu\|^2>\rho\qquad\forall\bu\in\cN. \ee Notice that by
\eqref{eq:delta psi restringida a N} and \eqref{eq:away from zero},
\cite[Proposition 6.7]{am} proves that $\cN$ is a Natural constraint
of $\Phi$, i.e., $\bu\in \h\setminus\{\bo\}$ is a critical point of
$\Phi$ if and only if $\bu$ is a critical point of $\Phi$
constrained on $\cN$.
\begin{Remarks}\label{re:functionl constrain y PS condition}
\hspace{1cm}
\begin{enumerate}
\item[(i)] The functional constrained on $\cN$ takes the form
\begin{equation}\label{eq:Functional constrain N2}
\Phi|_{\cN}(\bu)= \tfrac 16\|\bu\|^2+\tfrac{1}{12}\int_{\R^N} u^4\,dx.
\end{equation}
Even more, using \eqref{eq:away from zero} and \eqref{eq:Functional constrain N2},
\begin{equation}\label{eq:Functional constrain N}
\Phi(\bu)>\tfrac{1}{6}\rho\qquad \forall\,\bu\in\cN.
\end{equation}
Therefore, $\Phi$ is bounded from below on $\cN$, so we can try to
minimize it on the Nehari manifold.
\item[(ii)] Let us define
\[
2^*=\left\lbrace
\begin{array}{ccl}
\frac{2N}{N-4} & \text{if} & N>4, \\
\infty & \text{if} & 1\le N\le 4.
\end{array}
 \right.
\]
One has the following Sobolev embedding
\[
E\hookrightarrow L^p(\R^N),\qquad \text{for
}\left\{\begin{array}{rcl}2\le
p\le  2^*,& &  \mbox{if } N\neq 4\\
2\le p<2^*,& & \mbox{if } N= 4,
\end{array}\right.
\]
see for instance, \cite{Lions-JFA82,Adfou}.

In particular, this embeddings show that the functional $\Phi$ is
well defined for every $1\leq N\leq 7$.

Concerning the Palais-Smale condition for $2\le N\le 7$, (see Lemma
\ref{Lemma PS}) we will use that if  $N\ge 2$, replacing $E$ by the
radial subspace $H$, we have the following compact embedding
\[
H\hookrightarrow\hookrightarrow L^p(\R^N),\qquad \text{for }2<p<
2^*.
\]
The one dimensional  case ($N=1$) is analyzed in a different manner
in Subsection \ref{sec:4.2} because of the lack of compactness.
\end{enumerate}
\end{Remarks}

\begin{Remark}\label{semi-trivial}
System  \eqref{eq:NLS-KdV2} only admits one kind of semi-trivial
solutions of the form $(0, v)$. Indeed, if we suppose  $v=0$, the
second equation in \eqref{eq:NLS-KdV2} gives us that $u=0$ as well. Thus, let
us take $\bv_2=(0, V_2)$, where $V_2$ can be taken as a positive
radially symmetric  ground state solution of the equation
$\Delta^2v+\lambda_2v = \frac{1}{2}|v|v$. In particular, we can
assume that $V_2$ is positive because in other case, taking $|V_2|$,
it has the same energy. Moreover, if we denote by $V$ a positive radially
symmetric ground state solution of the equation $\Delta^2v+v =
\frac{1}{2}|v|v$, then, after some rescaling $V_2$ can be defined by \be\label{elemento
reescalado} V_2(x)=\lambda_2V(\sqrt[4]{\lambda_2}x). \ee As a
consequence, $\bv_2=(0,V_2)$ is a non-negative semi-trivial solution
of \eqref{eq:NLS-KdV2}, independently of the value of $\beta$.
\end{Remark}
We define the Nehari manifold corresponding to the single second
equation of \eqref{eq:NLS-KdV2} by
\[
\cN_2=\left\lbrace
v\in H\setminus\{0\}:J_2(v)=0
\right\rbrace \]
where
$$J_2(u):=I'_2(u)[u].$$
Let us define the tangent space to $\cN$ on $\bv_2$ by
\[
T_{\bv_2}\cN:=\left\lbrace
\bh\in \E:\Psi'(\bv_2)[\bh]=0
\right\rbrace,
\]
equivalently we define the tangent space to $\cN_2$ on $V_2$ by
\[
T_{V_2}\cN_2:=\left\lbrace
h\in E:J'_2(V_2)[h]=0
\right\rbrace.
\]
We can see that the  following equivalence holds:
\begin{equation}\label{eq:tang1}
\bh=(h_1,h_2)\in T_{\bv_2} \cN  \Longleftrightarrow h_2\in T_{V_2} \cN_2,
\end{equation}
in fact,
\[
\begin{array}{ccl}
\bh\in T_{\bv_2} \cN & \Longleftrightarrow & \Psi'(\bv_2)[\bh]=0\\
&\Longleftrightarrow& 2\bra V_2,h_2\ket_2-\frac{3}{2}\int_{\R^N}V_2^2h_2=0\\
&\Longleftrightarrow& J_2'(V_2)[h_2]=0\\
& \Longleftrightarrow & h_2\in T_{V_2} \cN_2.\\
\end{array}
\]

If we denote by $D^2\Phi_{\cN}$ the second derivative of $\Phi$
constrained on $\cN$, using that $\bv_2$ is a critical point of
$\Phi$, plainly we obtain that \be
D^2\Phi_{\cN}(\bv_2)[\bh]^2=\Phi''(\bv_2)[\bh]^2\qquad \forall\
\bh\in T_{\bv_2}\cN. \ee In the following result we establish the
character of $\bv_2$ in terms of the size of the coupling parameter.
\begin{Proposition}\label{Prop:fund}
There exists $\L>0$ such that:
\begin{itemize}
\item[(i)] if $\beta<\L$, then $\bv_2$ is a strict local minimum of $ \Phi$ constrained on $\cN$.
\item[(ii)] if $\beta>\L$, then $\bv_2$ is a saddle point of $ \Phi$ constrained on $\cN$. Moreover,
\be\label{eq:infimo B>L}
\inf\limits_\cN\Phi<\Phi(\bv_2).
\ee
\end{itemize}
\end{Proposition}

\begin{pfn}
\begin{itemize}
\item[(i)] We define
\be\label{Lambda} \L:=\inf\limits_{\varphi\in
H\setminus\{0\}}\frac{\|\varphi\|^2_1}{\int_{\R^N}V_2\varphi^2\,dx}.
\ee For $\bh\in  T_{\bv_2}\cN$ one has that \be\label{eq:seg deriv}
D^2\Phi_\cN
(\bv_2)[\bh]^2=\Phi''(\bv_2)[\bh]^2=\|h_1\|^2_1+I_2''(V_2)[h_2]^2-\beta\int_{\R^N}V_2h_1^2\,dx.
\ee By \eqref{eq:tang1} $\bh=(h_1,h_2)\in
T_{\bv_2}\cN\Leftrightarrow h_2\in T_{V_2} \cN_2$. Then, using that
$V_2$ is a minimum of $I_2$ on $\cN_2$, there exists a constant
$c_2>0$ such that \be\label{eq:minimo-pos}
 I_2'' (V_2)[h_2]^2\ge c_2\|h_2\|_2^2.
\ee Since $\beta<\Lambda$, \eqref{Lambda} and \eqref{eq:seg deriv}
there exists $c_1>0$ such that \be\label{eq:desig seg der}
D^2\Phi_\cN (\bv_2)[\bh]^2\geq c_1\|h_1\|_1^2+c_2\|h_2\|_2^2, \ee
proving that $\bv_2$ is a strict local minimum of $\Phi$ on $\cN$.
\item[(ii)]
Since $\beta >\L$, there exists $\wt{h}\in H$ such that
$$
\L< \frac{\|\wt{h}\|_1^2}{\int_{\R^N} V_2\wt{h}^2\,dx}<\b.
$$
Then, taking $\bh_1=(\wt{h},0)\in T_{\bv_2}\cN$ it yields
$$
D^2\Phi_\cN (\bv_2)[\bh_1]^2=\|\wt{h}\|_1^2 -\b_{\R^N} V_2
\wt{h}^2dx<0,$$ and  taking $h_2\in T_{V_2}\cN_2$ not equal to zero,
then $\bh_2=(0,h_2)\in T_{\bv_2}\cN$ and
$$D^2\Phi_\cN
(\bv_2)[\bh_2]^2=I_2'' (V_2)[h_2]^2\ge c_2\|h_2\|_2^2>0.
$$
Therefore, $\bv_2$ is a saddle point of $\Phi$ on $\cN$ and obviously inequality \eqref{eq:infimo B>L} holds.
\end{itemize}
\end{pfn}

To conclude this section we also prove that the functional $\Phi$ satisfies the PS condition constrained to $\cN$ on the high-dimensional case.
\begin{Lemma}\label{Lemma PS}
Assume that $2\le N\leq 7$, then $\Phi$ satisfies the PS condition constrained on $\cN$.
\end{Lemma}
\begin{pfn}
Let $\bu_n=(u_n,v_n)\in \cN$ be a PS sequence, i. e., \be\label{PS}
\Phi(\bu_n)\to c \quad \text{and}\quad \nabla_\cN\Phi(\bu_n)\to
0,\quad \text{as}\quad n\to\infty. \ee
 From \eqref{eq:Functional constrain N2} and the first convergence in \eqref{PS} it follows that $\bu_n$ is bounded, then we
 have a weakly convergent subsequence (denoted equals for short) $\bu_n\rightharpoonup \bu_0\in \mathbb{H}$.
 Since $H$ is compactly embedding into $L^p(\R^N)$ for $2<p<4+\frac{2}{3}$ and $2\le N\leq 7$
 (see Remark \ref{re:functionl constrain y PS condition}-(ii)), we infer that
\[
\int_{\R^N}u_n^4\,dx\to\int_{\R^N}u_0^4\,dx,
\qquad \int_{\R^N}|v_n|^3\,dx\to\int_{\R^N}|v_0|^3\,dx,
\qquad \int_{\R^N}u_n^2v_n\,dx\to\int_{\R^N}u_0^2v_0\,dx.
\]
Moreover, using the fact that $\bu_n\in\cN$ and \eqref{eq:away from zero}, we have
\[
\|\bu_n\|^2=\intN u_n^4\,dx +\tfrac{1}{2} \intN |v_n|^3\,dx
+\tfrac{3}{2}\beta\intN u_n^2v_n\,dx\to\int_\R u_0^4
\,dx+\tfrac{1}{2} \intN |v_0|^3\,dx +\tfrac{3}{2}\beta\intN
u^2_0v_0\,dx\geq\rho,
\]
which implies that $\bu_0\neq \bo$. The constrained gradient
satisfies
\begin{equation}\label{four2}
\nabla_\cN\Phi(\bu_n)=\Phi'(\bu_n)-\l_n\Psi'(\bu_n)\to 0,
\end{equation}
then, taking into account \eqref{eq:delta psi restringida a N},
\eqref{eq:away from zero}, the fact that
$\Phi'(\bu_n)[\bu_n]=\Psi(\bu_n)=0$, and evaluating the identity of
expression \eqref{four2} at $\bu_n$ we deduce that $\l_n\to 0$ as
$n\to\infty$. We also have that $\|\Psi'(\bu_n)\|$ is bounded.
Hence, from \eqref{four2}, jointly with the fact $\l_n\to 0$, we
obtain
$$
\|\Phi'(\bu_n)\|\le\|\nabla_\cN \Phi(\bu_n)\|+|\l_n|
\|\Psi'(\bu_n)\|\to 0\qquad\text{as}\quad n\to\infty.
$$
To finish the proof, since $\Phi'(\bu_n)[\bu_0]\to 0$ as $n\to
\infty$, it follows that $\bu_n\to \bu_0$ strongly.
\end{pfn}

\section{Existence results}\label{sec:4}
This section is divided into two subsections depending on the dimension of problem \eqref{eq:NLS-KdV2}.
\subsection{High-dimensional case, $2\le N\le 7$.}\label{sec:4.1}

\

In this subsection we will see that the infimum of $\Phi$ constrained on the radial Nehari manifold, $\cN$, is attained under appropriate
parameter conditions.
We also prove the existence of a mountain pass critical point.

\begin{Theorem}\label{th:minimo1} Suppose $\beta>\L$ and $2\leq N\leq 7$. The infimum of $\Phi$ on $\cN$ is attained at some point $\wt{\bu}\geq \bo$
with $\Phi(\wt{\bu})<\Phi(\bv_2)$ and both components $\wt{u},\wt{v}\not \equiv 0$.
\end{Theorem}
\begin{pfn} By the Ekeland's Variational Principle (see \cite{eke} for further details) there exists a minimizing PS sequence $\bu_n\in\cN$, i.e.,
\[
\Phi(\bu_n)\to c:=\inf\limits_\cN\Phi\quad\text{and}\quad
\nabla_\cN\Phi(\bu_n)\to 0.
\]
Due to the Lemma \ref{Lemma PS}, there exists $\wt{\bu}\in\cN$ such that
$$
\bu_n\to \wt{\bu}\quad \text{strongly as}\quad n\to\infty,
$$
 hence $\wt{\bu}$ is a minimum point of $\Phi$ on $\cN$. Moreover, taking into account Proposition \ref{Prop:fund}-(ii), we have:
\[
\Phi(\wt{\bu})=c<\Phi(\bv_2).
\]
Note that the second component $\wt{v}$ can not be  zero, because if that occur then $\wt{\bu}\equiv 0$ due to the form of the second equation of \eqref{eq:NLS-KdV2}, and zero is not in $\cN$. On the other hand, if we suppose that the first component $\wt{u}\equiv 0$, then
$$
I_2(\wt{v})=\Phi(\wt{\bu})<\Phi(\bv_2)=I_2(V_2),
$$
and this is a contradiction with the fact that $V_2$ is a ground state of the equation $\Delta^2v+\l_1v = \frac{1}{2}|v|v$.

In general we can not ensure that both components of $\wt{\bu}$ are non-negative, thus, in order to obtain this fact we take $t|\wt{\bu}|\in \cN$, and we will show that
 $$
 \Phi(t|\wt{\bu}|)\leq\Phi(\wt{\bu}).
 $$
 Note that by \eqref{eq:Functional constrain N} we have that
\begin{equation}\label{eq:comparar t}
\Phi(t|\wt{\bu}|)= \tfrac{1}{ 6}t^2\|\wt{\bu}\|^2+\tfrac{1}{12}t^4\int_{\R^N} \wt{u}^4\,dx,
\qquad\Phi(\wt{\bu})= \tfrac 16\|\wt{\bu}\|^2+\tfrac{1}{12}\int_{\R^N} \wt{u}^4\,dx.
\end{equation}
Hence, to prove $\Phi(t|\wt{\bu}|)\leq\Phi(\wt{\bu}) $ is equivalent to show that $t\leq 1$. Taking into account that $\Psi(t|\wt{\bu}|)=0$, we find:
\begin{align*}
0=\Psi(t|\wt{\bu}|)&=t^2\|\wt{\bu}\|^2-t^4\int_{\R^N} \wt{u}^4\,dx -\tfrac{1}{2}t^3 \int_{\R^N} |\wt{v}|^3\,dx -\tfrac{3}{2}t^3\beta\int_{\R^N} \wt{u}^2|\wt{v}|\,dx,
\end{align*}
which is equivalent to
\begin{equation}\label{eq:condicion con t}
0=\|\wt{\bu}\|^2-t^2\int_{\R^N} \wt{u}^4\,dx -\tfrac{1}{2}t \int_{\R^N} |\wt{v}|^3 \,dx-\tfrac{3}{2}t\beta\int_{\R^N} \wt{u}^2|\wt{v}|\,dx.
\end{equation}
Furthermore, since $\wt{\bu}\in \cN$ we also have,
\begin{equation}\label{eq:condicion sin t}
0=\Psi(\wt{\bu})=\|\wt{\bu}\|^2-\int_{\R^N} \wt{u}^4\,dx -\tfrac{1}{2} \int_{\R^N} |\wt{v}|^3 \,dx-\tfrac{3}{2}\beta\int_{\R^N} \wt{u}^2\wt{v} \,dx.
\end{equation}
Now, if we suppose that $t>1$ it follows that
\[
t^2\int_{\R^N} \wt{u}^4\,dx +\tfrac{1}{2}t \int_{\R^N} |\wt{v}|^3 \,dx+\tfrac{3}{2}t\beta\int_{\R^N} \wt{u}^2|\wt{v}|\,dx>\int_{\R^N} \wt{u}^4\,dx +\tfrac{1}{2} \int_{\R^N} |\wt{v}|^3\,dx +\tfrac{3}{2}\beta\int_{\R^N} \wt{u}^2|\wt{v}|\,dx.
\]
Then, thanks to \eqref{eq:condicion con t} we obtain
\begin{equation}\label{desig}
0<\|\wt{\bu}\|^2-\int_{\R^N} \wt{u}^4\,dx -\tfrac{1}{2} \int_{\R^N} |\wt{v}|^3\,dx -\tfrac{3}{2}\beta\int_{\R^N} \wt{u}^2|\wt{v}|\,dx.
\end{equation}
Combining \eqref{eq:condicion sin t} with \eqref{desig} we arrive at
\[0<\tfrac{3}{2}\beta\int_{\R^N} \wt{u}^2\left(\wt{v}-|\wt{v}|\right)\,dx,\]
which is a contradiction. Consequently, $t\leq 1$ and therefore $\Phi(t|\wt{\bu}|)\leq\Phi(\wt{\bu})$. On the other hand, we know that $\Phi$ attains its
infimum at $\wt{\bu}$ on $\cN$ and, therefore, the last inequality can not be strict. Moreover, due to  \eqref{eq:comparar t} it can not happen that $t<1$ and, hence, $ t=1$ and
$$\Phi(|\wt{\bu}|)=\Phi(\wt{\bu}).$$
Redefining $\wt{\bu}$ as $|\wt{\bu}|$ we finally have that the minimum on the Nehari manifold is attained at $\wt{\bu}\geq 0$ with non-trivial components.
\end{pfn}

\begin{Theorem}\label{Th:lambda>L_2}
Assume  $2\leq N\leq 7$, $\beta>0$. There exists a positive constant
$\L_2$ such that, if $\l_2>\L_2$, the functional $\Phi$ attains its
infimum on $\cN$ at some $\widehat{\bu}\geq\bo$ with
$\Phi(\widehat{\bu})<\Phi(\bv_2)$ and both
$\widehat{u},\widehat{v}\not \equiv 0$.
\end{Theorem}

\begin{pfn}
Using the same argument as above in the  Theorem \ref{th:minimo1},
we prove that the infimum is attained at some point
$\widehat{\bu}\in \cN$, but to show that
$\widehat{u},\widehat{v}\not \equiv 0$ we need to ensure that
$\Phi(\widehat{\bu})<\Phi(\bv_2)$. In Theorem \ref{th:minimo1} this
fact was proved for the case $\beta>\L$ and here we need to prove it
for $0<\beta\leq\L$. In this case the point $\bv_2$ is a strict
local minima and this does not guarantee that $\widehat{\bu}\not
\equiv\bv_2$.

Then, to see $\Phi(\widehat{\bu})<\Phi(\bv_2)$ we will use a similar procedure to the one applied in \cite{c2} showing
that there exists an element of the form
$$
\bw=t(V_2,V_2)\in\cN\quad \hbox{with}\quad \Phi(\bw)<\Phi(\bv_2),
$$
for $\lambda_2$ big enough.

Notice that, thanks to the equation $\Psi(\bw)=0$ we have that any
$t>0$ satisfies the following condition \be\label{eq:condicion1}
t^2\|(V_2,V_2)\|^2-t^4\intN V_2^4\,dx
-\tfrac{1}{2}t^3(1+3\beta)\intN V_2^3\,dx=0, \ee and by definition
we also have \be\label{eq:norma doble}
\|(V_2,V_2)\|^2=2\|V_2\|_2^2+(\lambda_1-\lambda_2)\intN V_2^2\,dx.
\ee Moreover, since $V_2\in\cN_2$, we have \be\label{eq:norma
simple} \|V_2\|_2^2-\tfrac{1}{2}\intN V_2^3\,dx=0. \ee Substituting
\eqref{eq:norma doble}  and \eqref{eq:norma simple} in
\eqref{eq:condicion1} it follows \be\label{eq:condicion2}
t^2\left(\intN V_2^3\,dx+(\lambda_1-\lambda_2) \intN
V_2^2\,dx\right) -t^4\intN V_2^4\,dx -\tfrac{1}{2}t^3(1+3\beta)\intN
V_2^3\,dx=0. \ee Hence, applying the rescaling \eqref{elemento
reescalado}  yields \be\label{cambio p} \intN
V_2^p\,dx=\lambda_2^{p-\frac{N}{4}}\intN V^p\,dx. \ee Subsequently,
substituting \eqref{cambio p} for $p=2,3,4$ into
\eqref{eq:condicion2} and dividing by $t^2\lambda_2^{3-\frac{N}{4}}$
we have that \be\label{eq:condicion} \intN
V^3\,dx+\dfrac{\lambda_1-\lambda_2}{\lambda_2}\intN
V^2\,dx-t^2\lambda_2\intN V^4\,dx-\tfrac{1}{2}t(1+3\beta)\intN
V^3\,dx=0. \ee Moreover, due to \eqref{eq:Functional constrain N2},
\eqref{eq:norma doble} and \eqref{eq:norma simple} we find
respectively the expressions \be\label{eq: forma de Phi(w)}
\Phi(\bw)=\tfrac{1}{6}t^2\left( \intN
V_2^3\,dx+(\lambda_1-\lambda_2) \intN V_2^2\,dx\right)
+\tfrac{1}{12}t^4\intN V_2^4\,dx, \ee \be\label{eq: forma de
Phi(bv_2)}
\Phi(\bv_2)=I_2(V_2)=\tfrac{1}{2}\|V_2\|_2^2-\tfrac{1}{6}\intN
V_2^3=\tfrac{1}{12}\intN V_2^3. \ee Furthermore, we are looking for
the inequality $\Phi(\bw)< \Phi(\bv_2)$, or equivalently,
\be\label{40} \tfrac{1}{6}t^2\left( \intN
V_2^3\,dx+(\lambda_1-\lambda_2) \intN V_2^2\,dx\right)
+\tfrac{1}{12}t^4\intN V_2^4\,dx- \tfrac{1}{12}\intN V_2^3\,dx <0,
\ee and then, applying again \eqref{cambio p} and multiplying
\eqref{40} by $6\lambda^{\frac{N}{4}-3}$, we actually have
\be\label{desigualdad de w} t^2\left( \intN
V^3\,dx+\dfrac{\lambda_1-\lambda_2}{\lambda_2}\intN V^2\,dx\right)
+\tfrac{1}{2}t^4\lambda_2\intN V^4\,dx-\tfrac{1}{2}\intN V^3\,dx<0.
\ee
Solving \eqref{eq:condicion} the corresponding will provide us \eqref{desigualdad de w} for $\l_2$ large enough.

Therefore, there exists a positive constant $\L_2$ such that for $\l_2>\L_2$ inequality \eqref{desigualdad de w} holds, proving that
$$
\Phi(\widehat{\bu})\leq\Phi(\bw)< \Phi(\bv_2).
$$
Finally, to show that $\widehat{\bu}\geq\bo$ and
$\widehat{u},\widehat{v}\not \equiv 0$ we can use the same argument
as in Theorem \ref{th:minimo1}.
\end{pfn}

In the following we will prove the existence of a MP critical point of $\Phi$ on $\cN$.
\begin{Theorem}\label{Montain Pass}
Assume  $2\leq N\leq 7$ and $\beta<\L$. There exists a constant
$\L_2$ such that, if $\l_2>\L_2$, then $\Phi$ constrained on $\cN$
has a Mountain-Pass critical point $\bu^*$ with
$\Phi(\bu^*)>\Phi(\bv_2)$.
\end{Theorem}
\begin{pfn}
Due to Proposition \ref{Prop:fund}-(i), $\bv_2$ is a strict local minima of $\Phi$ on $\cN$, and taking into account Theorem \ref{Th:lambda>L_2} we
obtain $\L_2$ such that, for $\lambda_2>\L$, we have $\Phi(\widehat{\bu})<\Phi(\bv_2)$. Under those conditions
we are able to apply the Mountain Pass Theorem (see \cite{ar} for further details) to $\Phi$ on $\cN$, that provide us with a PS sequence $\bv_n\in\cN$ such that
\[
\Phi(\bv_n)\to m:=\inf\limits_{\gamma\in\Gamma}\max\limits_{0\le t\le 1}\Phi(\gamma(t)),
\]
where
\[
\Gamma:=\left\{\gamma:[0,1]\to\cN\ \ \text{continuous}\  |\ \gamma(0)=\bv_2,\ \gamma(1)=\widehat{\bu} \right\}.
\]
Furthermore, applying the Lemma \ref{Lemma PS},
we are able to find a subsequence of $\bv_n$ such that (relabelling) $\bv_n\to\bu^*$ strongly in $\h$. Thus, $\bu^*$ is a critical point of $\Phi$
satisfying
\[
\Phi(\bu^*)>\Phi(\bv_2),
\]
which conclude the proof.
\end{pfn}
\subsection{One-dimensional case, $N=1$.}\label{sec:4.2}

\

Here we must  point out that we do not have the compact embedding even for $\h$. However, we  will show that for a PS sequence we are able to find a subsequence
for which its weak limit is a solution of \eqref{eq:NLS-KdV2} belonging to $\E$. Thus, in order to avoid the lack of compactness for $N=1$ we will use the following result of measure theory that one can find in \cite{lions2};
see also \cite{c-fract,c3} for an application of this procedure to a similar problem.
\begin{Lemma}\label{lem:measure}
If $2<q<\infty$, there exists a constant $C>0$ so that
\begin{equation}\label{eq:measure}
\intR |u|^q \, dx\le C\left( \sup_{z\in\R}\int_{|x-z|<1}|u(x)|^2\,dx\right)^{\frac{q-2}{2}}
\| u\|^2_{E},\quad \forall\: u\in E.
\end{equation}
\end{Lemma}

The next result is analogous to Theorem \ref{th:minimo1} for the
one-dimensional case and working on the full Nehari manifold
$\mathcal{M}$ defined by \eqref{gemanif}.

\begin{Theorem}\label{th:minimo N=1}
Suppose $N=1$ and $\beta>\L$. The infimum of $\Phi$ on $\mathcal{M}$
is attained at some $\wt{\bu}\geq \bo$ with both components
$\wt{u},\wt{v}\not \equiv 0$. Moreover,
$\Phi(\wt{\bu})<\Phi(\bv_2)$.
\end{Theorem}
\begin{pfn}
 Again, by the Ekeland's variational principle there exists a PS sequence $\bu_n\in\mathcal{M}$, i.e.,
\[
\Phi(\bu_n)\to c:=\inf\limits_{\mathcal{M}}\Phi\quad\text{and}\quad
\nabla_{\mathcal{M}}\Phi(\bu_n)\to 0,
\]
such that, $\bu_n$ is bounded since \eqref{eq:Functional constrain
N2}. Also, we can assume that the sequence $\bu_n$ possesses a
subsequence such that (relabelling) it weakly converges
$\bu_n\rightharpoonup \bu$ in $\E$, $\bu_n\to \bu$ strongly in
$\mathbb{L}^q_{loc}(\R)=L^q_{loc}(\R)\times  L^q_{loc}(\R)$ for
every $1\le q<\infty$ and $\bu_k\to \bu$ a.e. in $\mathbb{R}$.
Moreover, arguing in the same way as in Lemma \ref{Lemma PS} we
obtain $\Phi'(\bu_n)\to 0$ as $n\to\infty$.

Furthermore, using the idea performed in \cite{c2} we will prove that there is no loss of mass at infinity for $\mu_n(x):=u_n^2(x)+v_n^2(x)$, where
 $\bu_n=(u_n,v_n)$, i.e, there exist $R, C>0$ such that
\be\label{eq:vanishing}
\sup_{z\in\R}\int_{|z-x|<R}\mu_n (x)\,dx\ge C>0,\quad\forall n\in\mathbb{N}.
\ee
On the contrary, if we suppose
$$
\sup_{z\in\R}\int_{|z-x|<R}\mu_k(x)\,dx\to 0,
$$
and thanks to Lemma \ref{lem:measure} applied in a similar way as in \cite{c-fract},  we find that $\bu_k\to \bo$ strongly in
$\mathbb{L}^{q}(\R)$ for any $2<q<\infty$, This is a contradiction since $\bu_n \in\cN$, and due to \eqref{eq:Functional constrain N} jointly with the fact $\Phi(\bu_n)\to c$ we have
$$
0< \frac 17\rho <c+o_n(1)=\Phi(\bu_n),\quad\mbox{with } o_n(1)\to 0
\quad\mbox{as }n\to\infty,
$$
hence \eqref{eq:vanishing} is true and there is no loss of mass at infinity.

We observe that  there is a sequence of points
$\{z_n\}\subset\R$ such that by \eqref{eq:vanishing}, the translated sequence $\overline{\mu}_n(x)= \mu_n(x+z_n)$ satisfies
$$
\liminf_{n\to\infty}\int_{B_R(0)}\overline{\mu}_n\,dx\ge C >0.
$$
Taking into account that $\overline{\mu}_n\to \overline{\mu}$
strongly in $L_{loc}^1(\R)$, we obtain that
$\overline{\mu}\not\equiv 0$, thus, the weak limit of
$\overline{\bu}_n(x):=\bu_n(x+z_n)$, which we denote it by
$\overline{\bu}$, is non-trivial. Notice that
$\overline{\bu}_n,\overline{\bu}\in\mathcal{M}$ and
$\overline{\bu}_n$ is $PS$ sequence of level $c$ for  $\Phi$ on
$\mathcal{M}$. Moreover, if we set $F=\Phi|_{\mathcal{M}}$
(similarly to \eqref{eq:Functional constrain N2}) and using Fatou's
lemma we obtain the following
$$
\Phi (\overline{\bu})  =  \dyle F(\overline{\bu})
 \le  \dyle\liminf_{n\to\infty} F(\overline{\bu}_n)
 =  \dyle\liminf_{n\to\infty}\Phi(\overline{\bu}_n)=
   \dyle\liminf_{n\to\infty}\Phi(\bu_n)= c.
$$
Therefore, $\overline{\bu}$ is a non-trivial critical point of $\Phi$ constrained on $\mathcal{M}$. Furthermore, it is not a semi-trivial solution
because of $\Phi(\overline{\bu})<\Phi(\bv_2)$ from Proposition \ref{Prop:fund}-$(ii)$. Finally, to show that $\overline{\bu}\geq\bo$
and both components $\overline{u},\overline{v}\not \equiv 0$, we apply the same argument used in Theorem \ref{th:minimo1}.
\end{pfn}

Theorem \ref{Th:lambda>L_2} can by extended to the one-dimensional
case directly using the same idea as we have performed in the last
proof, obtaining the following.
\begin{Corollary}\label{Cor:lambda>L_2}
Assume  $N=1$, $\beta>0$. There exists a positive constant $\L_2$
such that, if $\l_2>\L_2$, the functional $\Phi$ attains its infimum
on $\cN$ at some $\widehat{\bu}\geq\bo$ with
$\Phi(\widehat{\bu})<\Phi(\bv_2)$ and both
$\widehat{u},\widehat{v}\not \equiv 0$.
\end{Corollary}
To finish, for $N=1$, Theorem \ref{Montain Pass}  can be obtained in
a similar manner, obtaining the following.
\begin{Corollary}\label{Cor:Montain Pass}
Assume  $N=1$ and $\beta<\L$. There exists a constant $\L_2$ such
that, if $\l_2>\L_2$, then $\Phi$ constrained on $\cN$ has a
Mountain-Pass critical point $\bu^*$ with $\Phi(\bu^*)>\Phi(\bv_2)$.
\end{Corollary}

%%%%%%%%%%%%%%%%%%%%%%%%%%%%%%%%%%%%%%%%%%%%%%%%%%%%%%%%%%%%%%%%%%%%%
\noindent{\bf Acknowledgements}.
First author was partially supported by the Ministry of Economy and Competitiveness of
Spain under research project RYC-2014-15284. Second author was partially supported by Ministry of Economy and Competitiveness of Spain and FEDER funds,
under research project MTM2013-44123-P.
%%%%%%%%%%%%%%%%%%%%%%%%%%%%%%%%%%%%%%%%%%%%%%%%%%%%%%%%%%%%%%%%%%%%%

%%%%%%%%%%%%%%%%%%%%%%%%%%%%%%%%%

\end{document}